\documentclass[12pt,a4paper]{amsart}
\usepackage{amsmath}
\usepackage[left=3cm,right=3.5cm,top=4cm,bottom=4cm]{geometry} 
\usepackage{meintex}

\usepackage[algo2e,norelsize]{algorithm2e}

\usepackage{mdframed}
\usepackage[section]{placeins}
\usepackage{tikz}
\usepackage[font=small,labelfont=bf]{caption}

\begin{document}

\title[Space mapping for Stochastic Interacting Particle Systems]{Space mapping-based Receding Horizon Control for Stochastic Interacting Particle Systems: dogs herding sheep}

%

\author[R.~Pinnau]{Ren\'{e} Pinnau}
\address[R.~Pinnau]{TU Kaiserslautern \\ Erwin-Schr\"odinger-Str.\ 48, 67663 Kaiserslautern}
\email[]{pinnau@mathematik.uni-kl.de}

\author[C.~Totzeck]{Claudia Totzeck}
\address[C.~Totzeck]{TU Kaiserslautern \\ Erwin-Schr\"odinger-Str.\ 48, 67663 Kaiserslautern}
\email[]{totzeck@mathematik.uni-kl.de}

\maketitle

\begin{abstract}
Control of stochastic interacting particle systems is a non-trivial task due to the high dimensionality of the problem and the lack of fast algorithms. Here, we propose a space mapping-based approximation of the stochastic control problem by  solutions of the deterministic one. In combination with the receding horizon control technique this yields a reliable and fast numerical scheme for the closed loop control of stochastic interacting particle systems. As a numerical example we consider  the  herding of sheep with dogs. The  numerical results underline the feasibility of our approach and further show stabilizing behaviour of the closed loop control.
\end{abstract}

\section{Introduction}
\label{intro}
Collective behaviour of crowds or swarms has been investigated by various researchers in the past decades \cite{ChristianiPiccoliTosin,CouzinKrause,CouzinKrauseJames,CuckerSmale,ParrishHamner}. First, the focus was on the  simulation of large groups, like flocks of birds and schools of fish, and their attractive and repulsive self-interaction  \cite{CarrilloMartinPanferov,Dorsogna}. The resulting models are able to reflect major properties of the interaction such as flocking and the formation of mills \cite{CarrilloKlarRoth}. Further, the stability of these patterns was analysed \cite{AlbiBalagueCarrillovonBrecht,CarrilloHuangMartin,Dorsogna}. Later, the models were refined to take into account view cones or topographical aspects like walls \cite{viewcones1,walls,Haskovec}. To include a  random disturbance of the individuals' behaviour one  introduces an additive Brownian motion in the velocity component of the dynamics \cite{jose12,jabin17}. Mathematically, this changes the model from ordinary differential equations (ODEs) to a system of stochastic differential equations (SDEs). 

Based on this knowledge, the investigation of the interaction of crowds and external agents became of interest \cite{AlbiPareschi,escobedo16}. In particular, the idea of controlling crowds with the help of the external agents \cite{Schafe2,paperMPC}. The corresponding  optimal control problem (OCP) is then constrained by the dynamics of the respective ODE or SDE system. For the deterministic problem one can employ  standard techniques from variational calculus to derive the gradient of the cost functional and to implement a tailored iterative scheme to compute the controls.

 Unfortunately, these classical methods cannot be directly adapted to the stochastic problems \cite{pfeiffer18}. In fact, the stochastic influence forces the decoupled forward and backward equations of an deterministic optimal control problem, to be a fully coupled Forward-Backward SDE system involving a ghost process to capture the uncertain terminal condition,  see, e.g., \cite{Perkowski} for the derivation of such a system based on a Hamiltonian formulation. First, steps towards a numerical realization in special cases can be found in \cite{gong2017}. 

Here, we are interested in controlling the crowd over a large time horizon, such that open loop control is not appropriate. Instead we use the closed loop receding horizon control to allow for feedback during the time evolution (see also \cite{albi18}). To deal with the stochastic nature of the model we employ the space mapping approach \cite{ReviewEngineer,koziel15}, which allows for the control of a high fidelity model (here the stochastic one) by the optimization of a surrogate model (here the deterministic one). 

 The space mapping approach first came up in the engineering community \cite{Bandler94} as a tool to solve large scale optimization problems with the help of an easier surrogate model. Through the years the technique became well-established in engineering and has been also recognized by the  mathematics community for various  applications like radiative heat transfer, control the dispersion of particles in a fluid, dynamic compressor optimization of gas networks and optimal inflow control of transmission lines, see, e.g., \cite{hemker05,Goettlich1,HertyBanda,hintermueller05,Goettlich2,MarheinekePinnauResendiz,PinnauThoemmes}.
 
Here,  we consider as new application the herding of a crowd of sheep using dogs with repulsive influence on the crowd. The combination of the space mapping technique with the receding horizon control will finally allow for the construction of a tailored closed loop algorithm to control  interacting stochastic particle systems.

The manuscript is organized as follows: in the next section the details of a general class optimal control problems with SDE constraints are given. Then, the space mapping approach is discussed in Section \ref{sec:SpaceMapping} and the Aggressive Monte Carlo Space Mapping Algorithm is presented. We derive the first-order optimality system of the deterministic ODE model and the gradient of the respective reduced cost functional which is needed for the numerical implementation in Section~\ref{sec:OptCon}. The algorithms for the numerical investigation are described in Section~\ref{sec:Numerics}. We present a projected gradient method for the deterministic optimization and discuss the receding horizon procedure for the closed loop control of the stochastic particle system. 
The feasibility of our approach is underlined by the  numerical results presented in Section 6. We discuss a space-mapping approach based on a mean-field approximation in Section 7, before we give conclusions and an outlook in Section 8. 

\section{The Control Problem}\label{sec:state}
In this section we define the general class of control problems constrained by a stochastic interacting particle system.

\subsection{Stochastic Interacting Particle System}
Let $D$ denote the space dimension and $T$ the length of the time interval under consideration. The positions and velocities of the particles are represented by
$$ X_i \colon [0,T] \rightarrow \R^{D}, \qquad  V_i \colon [0,T] \rightarrow \R^{D}, \qquad i=1,\ldots,N, $$ and combined in the vectors 
\begin{align*}
X(t) &= (X_1(t), \dots, X_N(t))^T \in \R^{ND}, \\ V(t) &= (V_1(t), \dots, V_N(t))^T \in \R^{ND} 
\end{align*} 
for each $t \in [0,T]$, respectively. In analogy, we consider $M$ external agents having positions 
$$ a_m \colon [0,T] \rightarrow \R^D, \quad m = 1,\ldots,M,$$
with $ a(t) = (a_1(t), \dots, a_M(t))^T \in \R^{MD}$  for each $t \in [0,T]$.

 Their velocities $$ u_m \colon [0,T] \rightarrow \R^D, \quad m=1,\ldots,M,$$
 are  combined in $ u(t) = (u_1(t), \dots, u_M(t))^T$ for each $t \in [0,T]$. Later, they act as control functions. We assume $u \in L^2([0,T], \R^{MD})$ .
 
The self-organisation of the crowd and the interaction of the particles with the agents is modelled with the help of radially symmetric interaction potentials $$\Phi_1, \Phi_2 \colon \R_0^+ \to \R, \quad \Phi_j(|x|) = \Phi_j(r), \quad j=1,2.$$  
For the sake of well-posedness, we assume that their first and second derivative $$\nabla_r \Phi_j (r) =: G_j(r) \text{ and }\nabla_r^2 \Phi_j(r) =: H_j(r)$$ are locally Lipschitz and globally bounded, i.e.,~$\Phi_j \in \mathcal C^2_b(\R_0^+)$ for $j=1,2$. 

Further, we include a friction term with parameter $\alpha>0$ and additive stochastic noise with strength $\sigma\ge 0$ influencing the velocities of the individuals. The friction models the lethargy of the individuals, while  the stochasticity allows for disturbances of the surroundings, that are not considered explicitly. Let $B_t^i$, $i=1,\dots,N$, denote independent $D$-dimensional Brownian motions. Then, the stochastic state system is given by
\begin{subequations}\label{eq:state}
	\begin{align}
	\dd X_i &= V_i \dd t, \qquad i =1,\dots, n,\\
	\dd V_i &= \left( - \frac{1}{N} \sum_{j=1}^N G_1(|X_i - X_j|) - \sum_{m=1}^M G_2(|X_i -a_m|) - \alpha v_i\right) \dd t + \sigma \dd B_t^i, \\
	\dd a_m &= u_m \dd t, \qquad m = 1,\dots,M,
	\end{align}
\end{subequations}
supplemented with initial data $Y_0 := (X_0,V_0,a_0)$.  The full state is a random variable $Y = (X,V,a)$. 

\begin{rem}
Clearly, for $\sigma =0$ the above system reduces to an ODE system, which we are going to use as the surrogate model for the space mapping procedure. 

It would be interesting to generalize the approach for common noise situations, i.e., $B_t^i = B_t$ in order to model effects on the system as a whole, instead of single particles. This will have impacts on the cost-functional and also on the mean-field equation discussed later on. For simplicity, we restrict ourselves here to the case of individual noise. 
\end{rem}

\subsection{Well-posedness of the State Systems}\label{sec:well-posedness-state}
Assuming a maximal velocity $u_\text{max}$ for the agents, we define the set of admissible controls 
$$ \Uad = \left\{ u \in L^2([0,T], \R^{MD}) \colon |u_m (t)| \le u_\text{max} \text{ for a.e. }t,\; m=1,\dots,M \right\}.$$
Note, that the ODE for the agents can be solved explicitly for given  $u \in \Uad$ which yields
$$ a(t) = a_0 + \int_0^t u(s) \dd s.$$
Indeed, we get  an absolutely continuous function $a$, which can be plugged into the SDE system governing the dynamic of the crowd. Using the assumption $\Phi_j \in \mathcal C_b^2(\R_0^+),$ we obtain weak solutions of the stochastic system in the sense of It\^o due to standard  SDE theory, see, e.g., \cite{Oksendal}. Further, the state fulfills $Y \in \mathcal{C}([0,T],\R^{ND}).$

\begin{rem}
	Note, that the stochastic system can be generalized to space and time dependent $\sigma = \sigma(x,t)$ without any effect on the well-posedness as long as the following conditions are satisfied
		\[
		|\sigma(t,x)| \le C(1 + |x|), \qquad |\sigma(t,x) - \sigma(t,y)| \le D|x-y|
		\]
		for some positive constants $C,\,D$  and $x,y \in \mathbb{R}^D$ .
\end{rem}	

\begin{rem}
	Note, that $u \in L^1([0,T,],\mathbb{R}^{MD})$ would be enough regularity to obtain the absolutely continuous function $a.$  But to identify the gradient for the numerical algorithm later on, we need a Hilbert-space structure. Thus, we choose the stronger assumption for $U_\text{ad}.$
\end{rem}	

In the case $\sigma = 0$, we obtain a deterministic ODE system which attains a unique solution by standard results from ODE theory.
This allows us to define the \textit{control-to-state} map $\G_c$ which assigns to each $u \in \Uad$ the unique solution $y$ of the ODE system. In analogy, we define $\G_f(u)=Y$ for the solution of the the SDE system.
For better readability we refer to states of the ODE system with lower-case letters and states corresponding to the SDE system with upper-case letters.

\subsection{The Cost Functional} In general, cost functionals involving empirical quantities, like expectation, variance or other kind of moments of the particle crowd are appropriate for the space mapping approach. 

In the following, we consider a specific cost functional that is based on the expected trajectory of the centre of mass of the crowd reflecting the aim of our application, i.e., steering the crowd to a predefined destination $\Edes$. To do so we define a time dependent reference state $\bar Z \colon [0,T] \rightarrow \mathbb{R}^D.$ Similar to the approach in \cite{paperMPC}, we measure the spread of the crowd around $\bar Z$. In particular, due to the stochastic behaviour of the state system we use the expected paths $\E[X]$. 

This leads to  the following cost functional
\begin{equation}\label{eq:costfunctional}
J(Y,u ; \bar Z,\bar u) :=  \int_0^T \frac{1}{2N} \sum_{k=1}^N \| \E [X_k(t)] - \bar Z(t)\|^2  + \frac{\gamma}{2} \|u(t)- \bar u(t)\|_{\R^{MD}}^2 \dd t, 
\end{equation}
where the first term tracks the expected centre of mass of the crowd and penalizes its distance to the desired trajectory. The second term measures the control costs and is weighted with the parameter $\gamma >0$.

\begin{rem}
The predefined desired trajectory $\bar Z(t)$ and the reference velocities $\bar u$ are input parameters for the cost functional. In the space mapping procedure, $\bar Z$ shall be replaced by the expected centre of mass and $\bar u$ by the optimal control of the surrogate model.
\end{rem}

To sum up, the SDE constrained optimal control problem of consideration is given by \\[0.25cm]
\begin{mdframed}
	\begin{problem}\label{prob:sde}
		Find $u^* \in \Uad$ such that
		$$ u^* = \argmin_{u \in \Uad} J(Y, u ; \bar Z , \bar u) \text{ subject to } \eqref{eq:state} \text{ with initial condition } Y(0) = Y_0. $$ 
	\end{problem}
\end{mdframed}

\begin{rem}
The existence of an optimal control can be shown with standard techniques from variational calculus \cite{Pinnau}. In fact, an existence result  can be obtained for all sequentially weak lower semicontinuous and coercive cost functionals $J(Y,u)$. Note that, in general, we cannot expect its uniqueness due to the non-convexity which is introduced by the nonlinearity in the state system.
\end{rem}

\section{The Space Mapping Approach}\label{sec:SpaceMapping}

The direct solution of this SDE constrained optimal control problem is a non-trivial task. Nevertheless, we can exploit the fact that the deterministic ODE model is for small noise $\sigma$ a good approximation for the stochastic one in combination with the space mapping procedure.

The general idea of space mapping for optimization problems is to approximate a complex (fine) model by a simple (coarse) surrogate model such that its main features are still resolved and the coarse model allows for a fast optimization. In particular,  no gradient information of the fine model needs to be computed.  Space mapping goes back to \emph{Bandler} \cite{Bandler94} and an excellent introduction is given in the review \cite{ReviewEngineer} and the references therein.

Let $\mathcal G_f$ and $\mathcal G_c$ be two operators mapping the fine control and the coarse control to some observable, respectively. To get an approximation of the fine model optimization 
 
$$u_f^* = \argmin_{u \in \Uad} \left| \mathcal{G}_f(u)-\bar w \right|$$
for a desired value $\bar w$, one uses optimizers of the coarse model, i.e.,

$$ u_c^* = \argmin_{u \in \Uad} \left| \mathcal{G}_c(u)-\bar w \right|.$$

For a better  approximation the space mapping function
$$
\T: \Uad \to \Uad, \quad \T(u_f) = \argmin_{u\in \Uad} \left| \mathcal{G}_c(u) - \mathcal{G}_f(u_f) \right|
$$
is introduced, which assigns to an input $u_f$ 	of the fine model a control $u_c$ for the coarse model, yielding the best approximation of the fine model output $ \mathcal{G}_f(u_f)$ by the coarse model output $\mathcal{G}_c(u_c)$.

If the observables to the respective optimizers are similar, i.e., $\mathcal{G}_f(u_f^*) \approx \mathcal{G}_c(u_c^*)$, we expect that it holds
$$
\T(u_f^*) = \argmin_{u\in \Uad} \left| \mathcal{G}_c(u) - \mathcal{G}_f(u^*_f) \right| \approx \argmin_{u\in \Uad} \left| \mathcal{G}_c(u) - \bar w \right| = u_c^*.
$$
Indeed, the space mapping is fixed by the observable defined by the operators $\mathcal G_f$ and $\mathcal G_c.$ In the following we use
\[
\mathcal G_f(u) = J(Y,u ; \bar Z,\bar u) =  \frac{1}{2N} \sum_{i=1}^N \left \| \begin{pmatrix}   \E[X^i] - \bar Z \\ \sqrt{\gamma}(u - \bar u) \end{pmatrix} \right \|_{L^2(0,T)}^2
\]
and
\[
\mathcal G_c(u) = J(y,u ; \bar Z,\bar u) =  \frac{1}{2N} \sum_{i=1}^N \left \| \begin{pmatrix}   x^i - \bar Z \\ \sqrt{\gamma}(u - \bar u) \end{pmatrix} \right \|_{L^2(0,T)}^2
\]
and compute $\E[X]$ with the help of a Monte Carlo simulation as proposed in  \cite{MarheinekePinnauResendiz}. Therefore, it makes sense to set $\bar w = 0.$ That means, in the following, the stochastic interacting particle system $(\sigma>0)$ will act as the fine model, while the coarse model is given by the deterministic particle system $(\sigma=0)$.

Another possible choice could use the solution operators of the state equations, i.e., $\mathcal G_f = \G_f$ and $\mathcal G_c = \G_c$ with $\bar w$ being an desired state.

\begin{rem}
Note, that the space mapping function $\T$ might be formally set valued if the optimization problem admits multiple solutions. Assumptions on the models ensuring that $\T$ is well defined are discussed in detail in \cite{hemker05,hintermueller05}.
\end{rem}

In general, the space mapping function $\T$ is directly not accessible, such that  there are several approximations proposed in the literature \cite{Bandler94,hemker05,ReviewEngineer}. These update the controls of the fine models iteratively. For example, Aggressive Space Mapping (ASM) and Trust Region Aggressive Space Mapping (TRASM) borrow the idea from quasi-Newton methods to approximate the Jacobian with the help of Broyden-type matrices. On the other hand, Hybrid Aggressive Space Mapping (HASM) combines the classical space mapping method with classical optimization techniques (cf.\ \cite{ReviewEngineer}). 

We use the ASM approach for the numerical computations below. Hence, the update $h^k$ for the next iterate is given by
$$ B^k h^k = -(\T(u_f^k) - u_c^*), \qquad u_f^{k+1} = u_f^k + \rho\, h^h, $$
where $B^k$ is the $k$-th Broyden matrix iterate and $\rho>0$ the step-length.

For a smooth presentation of the algorithm, we define the expected centre of mass of the stochastic particle crowd as
\begin{equation}\label{eq:expectedCenterOfMass}
\bar X(t) = \E\left[ \frac1N \sum_{k=1}^N X^k(t)\right], \qquad t \in [0,T]. 
\end{equation}

The resulting Aggressive Monte Carlo Space Mapping (AMCSM) approach\cite{MarheinekePinnauResendiz} is stated  in all details in Algorithm~\ref{AMCSM}.
\RestyleAlgo{boxruled}
\begin{algorithm2e}[!ht]
	\caption{Aggressive Monte Carlo Space Mapping (AMCSM)}\label{AMCSM}
	\KwData{initial values and parameters}
	\KwResult{control $u_f^*$}
	initialize counter $k=0$, approximate Jacobian $B^0 = I$ and tolerance $\epsilon_\text{SM}$\;
	Compute $u_f^0 = u_c^* = \argmin_{u_c} J(Z,u_c; \bar Z,0)$ subject to the deterministic model\;
	\While{ $ \|   \T(u_f^*) - u_c^* \| / \| u_c^* \|> \epsilon_\text{SM}$ }{
		evaluate the expected center of mass $\bar X$ given in \eqref{eq:expectedCenterOfMass} using MC simulations\;
		perform coarse model optimization $$ T(u_f^k) = \argmin_{u_c} J(Y,u_c ; \bar X, u_c^*) $$
		\If{$k>1$}{
			compute $B^k = B^{k-1} + (( \T(u_f^k) - u_c^*) \otimes h^{k-1}) / \left|h^{k-1}\right|^2$\;
		}
		solve $B^k h^k = -(\T(u_f^*) - u_c^*)$ for the update $h^k$\;
		update the control $u_f^{k+1} = u_f^k + h^k$\;
	}
\end{algorithm2e}

\begin{rem}
For the present application of dogs herding sheep, we need just one solve of the fine stochastic model, which involves the expensive Monte Carlo simulation in each step of the algorithm. The optimization step is only involving the coarse deterministic model, for which fast numerical algorithms based on gradient information are available. Clearly, for $N\gg 1$, the determinsitic problem is still challenging, one idea is to use a mean-field appoximation in this case. See the discussion on the mean-field limit in Section \ref{sec:meanfield}.
Further note, that the space-mapping approach discussed here does in general not yield a \emph{perfect} space mapping, such that the algorithm might terminate with a suboptimal solution (c.f.~\cite{hemker05}). This does not matter in our case, since we are designing a closed loop control with the help of the receding horizon control technique. The numerical results below indicate that the method proposed here, works fine for the problem at hand. Nevertheless, for other problems with short time horizons the space mapping solutions may fail to be robust. A qualitative study of the approximation and the robustness are subject to future work.
\end{rem}

\section{Optimal Control of the Coarse Model}\label{sec:OptCon}
The core of the space mapping approach is the fast  optimization of the coarse model. Since we intend to use a steepest descent algorithm, we derive the first-order optimality conditions for the coarse optimization problem. The derived adjoint information can then be used for the evaluation of the gradient of the reduced cost functional.

\subsection{First-Order Optimality Condition}

 For the deterministic optimal control problem with ODE constraints we can derive the adjoint system and the optimality condition with the help of the extended Lagrangian. Note, that the calculations are very similar to \cite{paperMPC}. 
The deterministic system 
\begin{subequations}\label{eq:detstate}
	\begin{align*}
	\frac{\dd }{\dd t} x_i &= v_i, \qquad \qquad i=1,\dots, N,\\
	\frac{\dd}{\dd t} v_i &= -\left( \frac{1}{N} \sum_{k=1}^N G_1(x_i - x_k) + \sum_{m=1}^M G_2(x_i -d_m) + \alpha v_i \right) =: -W^i(y), \\
	\frac{\dd}{\dd t} d_m &= u_m, \qquad \qquad m = 1,\dots,M,
	\end{align*}
\end{subequations}
can be compactly denoted by $\frac{\dd}{\dd t} y = F(y,u)$ , supplemented with the initial conditions $y(0) = y_0$.

We define the set of controls $\U$ and the state space $\Y$ as
\begin{equation*} 
\U = \{ u \in L^2([0,T], \R^{MD})\} , \quad  \Y = [H^1( [0,T],\mathbb{R}^{ND})]^2 \times H^1( [0,T],\mathbb{R}^{MD}).
\end{equation*}
Obviously it holds $\Uad \subset \U.$ Further, we define $\X:=[L^2( [0,T],\mathbb{R}^{ND})]^2\times L^2( [0,T],\mathbb{R}^{MD})$ and
\[
\Z:=\X\times\big([\mathbb{R}^{ND}]^2\times \mathbb{R}^{MD}\big),
\]
as the space of Lagrange multipliers, with $\Z^*$ being its dual. 

We define the state operator $e\colon \Y \times \U \rightarrow \Z^*$ for deterministic ODE as
\begin{equation*}
e(y,u) = \begin{pmatrix} \frac{\dd}{\dd t}y - F(y,u) \\ y(0)-y_0 \end{pmatrix}
\end{equation*}
and the dual pairing
\begin{equation*}
\langle e(y,u),(\xi,\eta) \rangle_{\Z^*,\Z} = \int_0^T (\frac{\dd}{\dd t}y(t) - F(y(t),u(t))) \cdot \xi(t) \dd t + (y(0) - y_0) \cdot \eta.
\end{equation*}
Let $(\xi,\eta)\in \Z$ denote the Lagrange multiplier which is in fact the adjoint state. Then, the extended Lagrangian corresponding to the coarse problem reads
\begin{equation*}
\mathcal{L}(y,u,\xi,\eta;\bar Z,\bar u) = J(y,u;\bar Z, \bar u) + \langle e(y,u),(\xi,\eta) \rangle_{\Z^*,\Z}.
\end{equation*}
As usual the first-order optimality condition of the coarse problem is given by
\begin{equation*}
d \mathcal{L}(y,u,\xi,\eta;\bar Z, \bar u) =0. 
\end{equation*}
Following the standard approach from variational calculus for the derivation of the adjoint equations (cf.~\cite{Pinnau}), we obtain the following first order optimality system.

\begin{theorem}\label{ODEKKT}
	Let $(y^*,u^*)$ be an optimal pair. Then, the first-order optimality condition corresponding to the coarse problem reads
	\begin{equation}\label{eq:var_ode}
	\int_0^T \big( \gamma ( u^*(t) - \bar u(t) ) - \xi_3(t) \big) \cdot(u(t)-u^*(t)) \dd t\ge 0  \qquad \text{for all\; $u\in \mathcal{U}_{ad}$},
	\end{equation}
	where $\xi=(\xi_1,\xi_2,\xi_3)\in \Y$ satisfies the adjoint system given by
	\begin{equation}\label{eq:xi}
	\frac{\dd}{\dd t}\xi_1 = -d_{{\emph x}} { W}(y^*)[\xi_2]- \frac{1}{NT}(x - \bar Z), \quad \frac{\dd}{\dd t}\xi_2 = \xi_1 - \alpha \xi_2, \quad  \frac{\dd}{\dd t}\xi_3 = -d_{{\emph d}} { W}(y^*)[\xi_2],
	\end{equation}
	supplemented with the terminal conditions $\xi_1(T) = 0$, $\xi_2 (T) = 0$, $\xi_3(T) = 0$.
\end{theorem}

\begin{rem}
The variational inequality in \eqref{eq:var_ode} can be derived as well with the help of the Pontryagin maximum principle. In view of the numerical implementation, the inequality is not handy.  We therefore choose the Lagrangian approach here, leading to explicit expressions for the adjoint which can be used in the algorithm. Together with a projection onto the feasible set $U_{ad}$ we can design a projected  gradient-descent method for the optimization problem, see Algorithm~\ref{alg:OptCoarse}.	
\end{rem}

\subsection{Gradient of the Reduced Cost Functional} \label{sec:redGrad}
In this section we introduce the reduced cost functional for the coarse model constraint and formally calculate its gradient which we need for the descent algorithm.
Using  the control-to-state map $\G_c$ we define the reduced cost functional as
\begin{equation*}
\hat{J}(u) := J(\G_c(u),u;\bar Z , \bar u).
\end{equation*}
Assuming sufficient regularity for $\G_c$ we further derive the gradient of the reduced cost functional. Making use of the state equation $e(y,u) = 0$ we implicitly obtain $d\G_c(u)$  via
\begin{equation*}
0= d_u e(\G_c(u),u) = d_y e(\G_c(u),u)[d\G_c(u)] + d_u e(\G_c(u),u).
\end{equation*}
With the help of the adjoint equation
\begin{equation*}
(d_ye(y,u))^*[\xi] = - d_y J(y,u) 
\end{equation*}
we compute the G\^ateaux derivative of $\hat{J}$ in direction $h\in \U$
\begin{equation*}
d\hat{J}(u)[h] = \langle d_y J(y,u), d\mathcal{S}_c(u)[h]\rangle_{\Y^*,\Y} + \langle d_u J(y,u),h\rangle_\U = \langle \gamma (u - \bar u) -\xi_3, h \rangle_{\U}.
\end{equation*}
Since $\U$ is a Hilbert space, we may use the Riesz representation theorem to identify the gradient of the reduced cost functional as
\begin{equation}\label{eq:gradient}
\nabla \hat{J} (u) = \gamma (u - \bar u)-\xi_3.
\end{equation}
Now, we have all ingredients at hand to state the gradient descent method for the numerical simulations. 

\section{Numerical Schemes}\label{sec:Numerics}
The Aggressive Monte Carlo Space Mapping algorithm (AMCSM) proposed in Algorithm~\ref{AMCSM} uses solutions of the coarse optimal control problem and only evaluations of  the fine stochastic particle system. 

\subsection{Optimization Algorithm for the Coarse Model}
We solve the deterministic ODE systems of state and adjoint problem with the explicit Euler scheme.
In the optimal control loop for deterministic problem, we update the controls using nonlinear conjugate gradient (NCG) steps. 

The step size for the gradient update is obtained by a line search based on the Armijo  rule with projection (cf.\ \cite{Pinnau}).

These ingredients define the numerical scheme for the deterministic optimization stated in Algorithm~\ref{alg:OptCoarse}, where we denote by $u^{n}$ the control of the $n$-th optimization iteration. When the optimal solution of the coarse problem $u_c^*$ is found, we compute $\bar x= \frac{1}{N} \sum_{i=1}^N x_i$, where the  $x_i$ refer to the optimal positions extracted from $\G_c(u_c^*).$

\begin{algorithm2e}[!ht]
	\caption{Optimal Control Algorithm for the Coarse Problem}\label{alg:OptCoarse}
	\KwData{initial data for states and control, stopping tolerance $\epsilon_\text{opt}$, time steps $K$,  desired destination $Z^*$}
	\KwResult{optimal control $u$, optimal states $y$}
	initialize\;
	\While{ $\| u^{n+1} - u^n \|_{L^2}> \epsilon_\text{opt}$ }{
		solve deterministic state system \eqref{eq:detstate}\;
		solve adjoint problem given in \eqref{eq:xi}\;
		compute gradient corresponding to \eqref{eq:gradient}\;
		compute step size using the Armijo rule with projection\;
		update controls by nonlinear conjugate gradient\; 
	}
\end{algorithm2e}
	In our particular case, the projection $\mathcal{P}_\U$ has the explicit representation 
	\begin{equation}\label{eq:projection}
	\mathcal{P}_\U(h)(t) = \begin{cases}
	u_{\text{max}}\frac{h_m(t)}{|h_m(t)|} & \text{for\; $|h_m(t)|>u_{\text{max}}$}, \\
	h_m(t) & \text{otherwise},
	\end{cases}\quad 
	\end{equation}
for $m=1,\ldots,M$ and $t\in [0,T]$.

\subsection{Receding Horizon Control}

The appropriate time horizon for steering the crowd to the given destination depends on the distance of the crowd to the destination and might be large. 
Since the space mapping procedure is based on optimal controls, we need to store the full forward information to compute the adjoints. On large time intervals this leads to an extensive memory consumption. Additionally, having the application of dogs herding sheep in mind, an open loop approach is rather unrealistic. In reality, a dog will react on the current state of the crowd, such that it makes more sense to model the problem using a closed loop ansatz.

\begin{figure}[htb]
	\begin{center}
		\begin{tikzpicture}[ultra thick]
		\draw (0,0) -- (7.5,0);
		\draw [dotted]   (7.5,0) -- (8.5,0);
		\draw (8.5,0) -- (10.,0);
		
		\draw (0,-0.75) -- (0,0.75);
		\begin{scope}[thin]
		\draw (1.5,-0.5) -- (1.5,0.5);
		
		\end{scope}
		\begin{scope}[ultra thick]
		\draw (3,-0.5) -- (3,0.5);
		\draw[<->] (0,-0.8) -- (3,-0.8);
		\draw (1.5,-1.25) node {$I_1$};
		\draw[<->] (0,0.8) -- (1.5,0.8);
		\draw (.75,1.25) node {$u_*^1$};
		\draw[<->] (1.5,-1.8) -- (4.5,-1.8);
		\draw (3.,-2.25) node {$[\frac{I_1}{2}, \frac{I_2}{2}]$};
		\draw[<->] (1.5,0.8) -- (3,0.8);
		\draw (2.25,1.25) node {$u_*^2$};
		\draw[<->] (3,-0.8) -- (6,-0.8);
		\draw (4.5,-1.25) node {$I_2$};
		\end{scope}
		
		\begin{scope}[thin]
		\draw (4.5,-0.5) -- (4.5,0.5);
		
		\end{scope}
		\begin{scope}[ultra thick]
		\draw (6,-0.5) -- (6,0.5);
		\end{scope}
		\begin{scope}[thin]
		\draw (7,-0.5) -- (7,0.5);
		
		\draw (9,-0.5) -- (9,0.5);
		\end{scope}
		\draw (10,-0.75) -- (10,0.75);
		\end{tikzpicture}
	\end{center}
	\caption{Visualization of the receding horizon procedure. The first iteration computes the optimal control on the interval $I_1$. Only the first half of it, $u_*^1$, is accepted as optimal solution.  Then the optimal control on the interval $[I_1/2, I_2/2]$ is computed. The first half, $u_*^2$, is accepted and clued to $u_*^1$.  These two steps give us the optimal control on $I_1$. We proceed iteratively up to the terminal time $T$.}
	\label{receedinghorizon}
\end{figure}
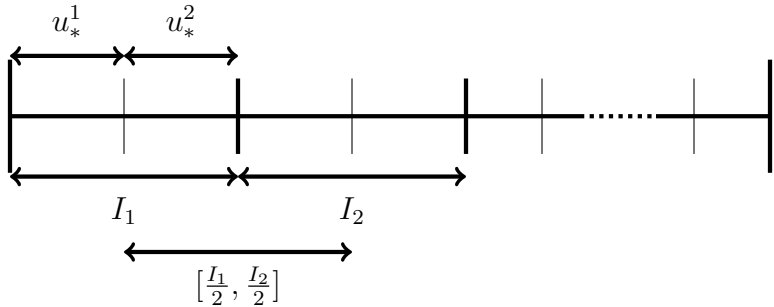

This is why a closed loop control for a large time horizon is preferable. Now, we are going to combine the above numerical approaches with the receding horizon control \cite{albi18}. In more detail, we split the time interval of interest $[0,T]$ into $K$ smaller intervals $I_1,\dots,I_K.$ Then, we apply the space mapping algorithm to these smaller intervals. In fact, we compute the stochastic output by an Euler-Mayurama scheme on $I_1$ but store only the first half of the solution. Then, we initialize the values using the optimal values at time $t = I_1/2$ and compute the solution on the interval $[\frac{I_1}{2}, \frac{I_2}{2}]$ and glue half of this solution to the one stored before. After two steps, we have the optimal control on the full interval $I_1$ available. We proceed iteratively until we reach the terminal time $T$. The receding horizon procedure is visualized in Figure~\ref{receedinghorizon}. Note that here is some freedom in choosing the length of the smaller interval. Numerical studies motivated us to use $\frac{I_k}{2}.$

\begin{rem}
Using this receding horizon procedure we need to adapt the desired trajectory $\bar Z$. Indeed, we cannot expect that the controls lead the crowd to the destination in one subinterval. Hence, we adapt $\bar Z$ on $I_k$ in the following way: we interpolate the distance of the initial centre of mass of the crowd and the desired destination $\Edes$ with the time steps used in one subinterval. Of course, this is not attainable for small $k$, nevertheless we simulate the deterministic optimal control using this interpolation as $\bar Z$ on $t \in [I_{k-1}, I_k]$. Then, we compute the trajectory of the center of mass corresponding to this solution. We expect this trajectory to be appropriate and use it in the space mapping procedure on the interval $[I_{k-1},I_k]$. 
\end{rem}

\section{Numerical Results}
In the following we present numerical results underlining the feasibility of our approach. In particular, we investigate the number of space mapping iterations needed to obtain appropriate results. Further, we shall see how the number of dogs is influencing the success of the herding procedure. Finally, we analyze numerically if the system is stabilized for large times $T\gg1$.

For the simulations we choose Morse potentials \cite{Dorsogna} to model the interaction:
\begin{gather*} G_j(|X_i - X_k|) = \nabla P_j(X_i,X_k), \; j\in \{1,2\}, \\  P_j(X_i,X_k) = C_{r,j} e^{-|X_i - X_k| / \ell_{r,j}} - C_{a,j} e^{-|X_i - X_k|/\ell_{a,j}}.  \end{gather*}
To realize the self-organization of the sheep we assume that they have some long range attraction and short range repulsion, i.e., we set $$ C_{r,1} = 1, C_{a,1} = 5e^{-4}, \ell_{r,1} = 2, \ell_{a,1}= 1e^{-2}.$$ Further, we assume the dogs to scare the sheep and therefore have stronger repulsive influence. This leads to $$ C_{a,2} = C_{a,1}, \ell_{a,2}=\ell_{a,1},  C_{r,2} = 1e^{-2}, \ell_{r,2} =0.5.$$

\begin{rem}
We emphasize that the space mapping control approach discussed here can be adapted to various other applications by changing the interaction potentials or the cost functional.
\end{rem}	
The following parameters are fixed for all simulations
$$ \gamma= 1e^{-2}, u_\text{max} = 5e^{-2}, K = 20, \epsilon_\text{opt}= 5e^{-3}, dt = 1e^{-2}, \alpha = 0.5,$$ where $dt$ denotes the time step size. Moreover, we choose $\sigma(x,t) = \sigma,$ i.e.~the stochastic force is independent of space and time. Nevertheless, $\sigma$ will be changed for different simulations and is thus specified explicitly later on as well as other parameters.\\

\subsection{Influence of the Stochasticity $\sigma$}

To study the influence of the stochasticity on the number of space mapping iterations, we set $$N=30, \quad M = 5,\quad  T = 20,$$ run $100$ Monte Carlo samples and  stop the iteration if
$$\| u_f - u_c^* \| / \| u_c^* \| < 0.3 \qquad \text{ or } \qquad  \| u_f^n - u_c^* \| / \| u_c^* \|  - \| u_f^{n+1} - u_c^* \| / \| u_c^* \| < 0.005$$ for two consecutive iterates $u_f^n$ and $u_f^{n+1}$. 

\begin{figure}[htb]
\centering
	\includegraphics[scale=0.5]{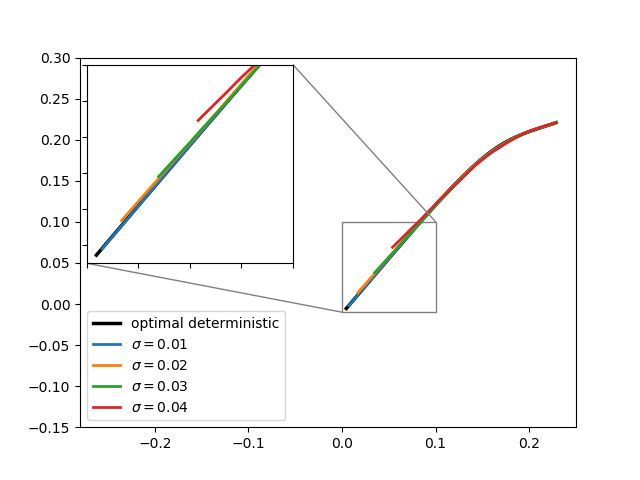}\\
	\includegraphics[scale=0.5]{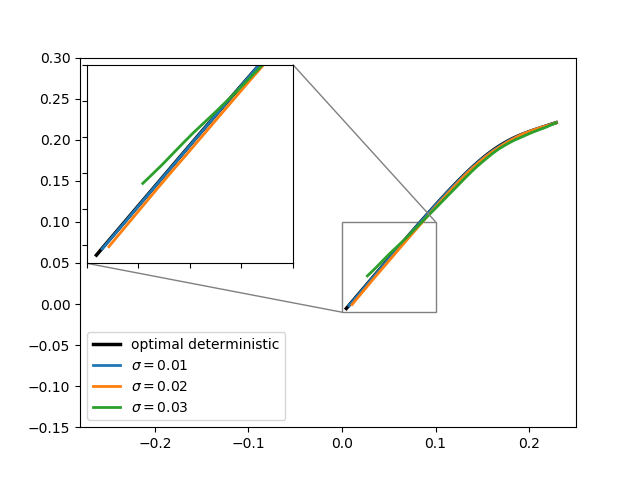}
	\caption{Up: We show the trajectories for the centre of mass of the crowd employing the optimal deterministic controls. The accuracy of the deterministic controls deteriorates as the stochastic influence increases. For $\sigma=0.03$ the stochastic influence starts to superimpose the crowd behaviour.  Down: The trajectories of the centre of mass of the crowd resuting from the space mapping procedure. We see that the trajectory corresponding to $\sigma = 0.02$ was improved.
	}
	\label{fig:compare_iterations}
\end{figure}

The accuracy of the deterministic controls deteriorates as the stochastic influence increases, see Figure~\ref{fig:compare_iterations} (up) as well as Table~\ref{tab:numerical_tests}. For $\sigma=0.03$ the stochastic influence starts to superimpose the crowd behaviour. Figure~\ref{fig:compare_iterations} (down) shows the trajectories of the center of mass of the crowd using space mapping. We see that space mapping works well for small values of $\sigma$. As the stochasticity starts to superimpose the crowd behaviour, the space mapping technique is not so efficient. This is  expected, since for large volatility the deterministic model is not a good approximation of the stochastic one. The second part of Table~\ref{tab:numerical_tests} shows results obtained with $1000$ Monte Carlo samples. The values change only slightly, which justifies to fix the number of MC samples to $100$ for the following computations.

\begin{rem}
We would like to emphasize that a basic Monte Carlo approach works fine in the present setting. For problems that are more involved, it may be necessary to use Multi-level Monte Carlo techniques in order to get efficient approximations for the stochastic states.
\end{rem}	

\begin{table}[htb]
	\begin{tabular}{l | c | c | c | c}
		MC = 100& $\sigma = 0.01$ & $\sigma = 0.02$ & $\sigma = 0.03$ & $\sigma = 0.04$ \\ \hline
		space mapping iterations & 0 & 1 & 3 & - \\
		$L^2$ error with deterministic control & $7.00\cdot e^{-3}$ & $3.73\cdot e^{-2}$ & $9.09\cdot e^{-2}$ & $1.57\cdot e^{-1}$ \\
		$L^2$ error after space mapping & $7.00\cdot e^{-3}$ & $1.05\cdot e^{-2}$ & $7.05 \cdot e^{-2}$ & - \\ \hline 
		MC = 1000& $\sigma = 0.01$ & $\sigma = 0.02$ & $\sigma = 0.03$ & $\sigma = 0.04$ \\ \hline
		space mapping iterations & 0 & 1 & 3 & - \\
		$L^2$ error with deterministic control & $9.28\cdot e^{-3} $ & $4.25\cdot e^{-2} $ & $9.70\cdot e^{-2} $ & $1.60\cdot e^{-1} $ \\
		$L^2$ error after space mapping & $9.28\cdot e^{-3}$ & $1.31\cdot e^{-2}$ & $7,07\cdot e^{-2}$ & - \\ \hline 
	\end{tabular}	
	\vspace{1em}
	\caption{Numerical investiation of the space mapping procedure. For $\sigma = 0.01$ no space mapping is needed, the optimal deterministic control is accepted. The number of space mapping steps increases with increasing stochastic strength. The $L^2$-error of the trajectory of the center of mass compared to the center of mass of the optimal deterministic solution increases as well for larger $\sigma$. The space mapping procedure is decreasing the error by a factor three for $\sigma = 0.02.$ As the stochastic starts to superimpose the crowd behaviour for $\sigma \ge 0.03$, we see that the space mapping approach decreases the error only marginally. The second part of the table shows results obtained with $1000$ Monte Carlo samples. The values change only slightly, which justifies to fix the number of MC samples to $100$ for the following computations. }
	\label{tab:numerical_tests}
\end{table}	

\subsection{Influence of the Number of Dogs $M$}
In the following figures, we depict sheep as blue dots, dogs as red triangles. The trajectories of the dogs are depicted as red lines and the trajectory of the center of mass of the crowd is the  blue line. A cross marks the desired location $Z_\text{des}.$ 

Varying the number of dogs leads to very different controls which can be visualized implicitly by the trajectories of the dogs. For this study we chose the parameter values
$$ \epsilon_\text{SM} = 0.5, \quad N=20, \quad  \sigma = 0.01. $$ 
Moreover, instead of fixing $T$ we used $| \bar X - Z_\text{des}| < 0.05$ as stopping criterion and did $100$ Monte Carlo runs. The change of the stopping criterion is necessary because we expect that a different number of dogs will need different times to steer the crowd to the desired destination.

\begin{figure}[htb]
	\includegraphics[scale=0.38]{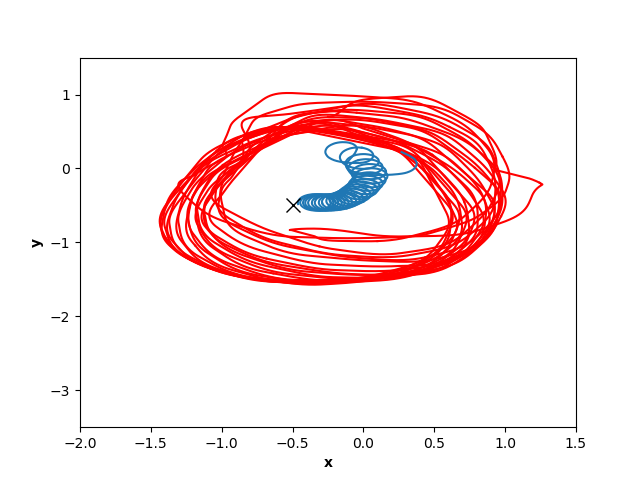} 
	\includegraphics[scale=0.38]{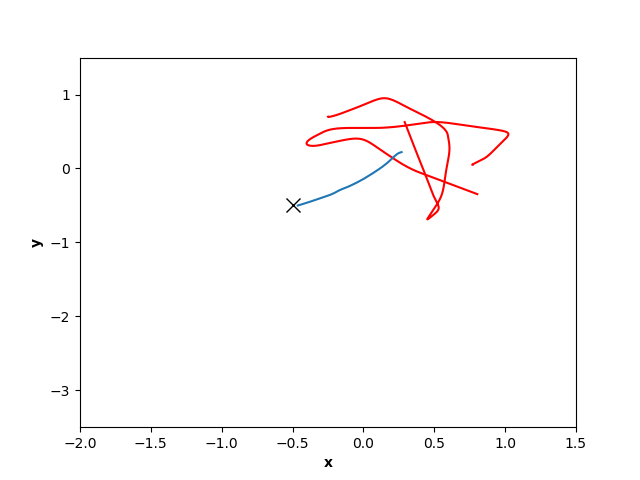}\\
	\includegraphics[scale=0.38]{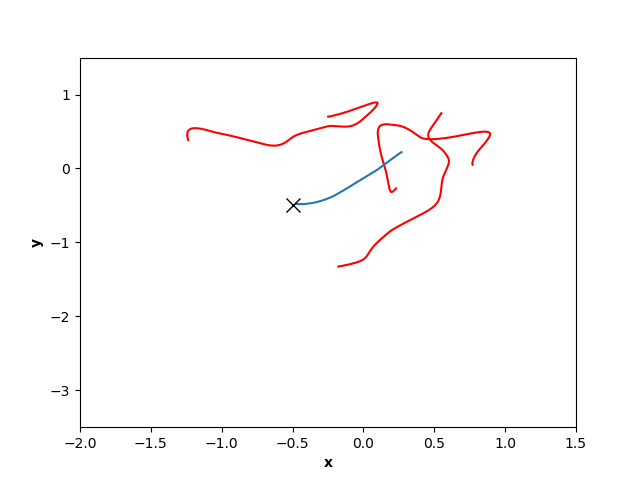}
	\includegraphics[scale=0.38]{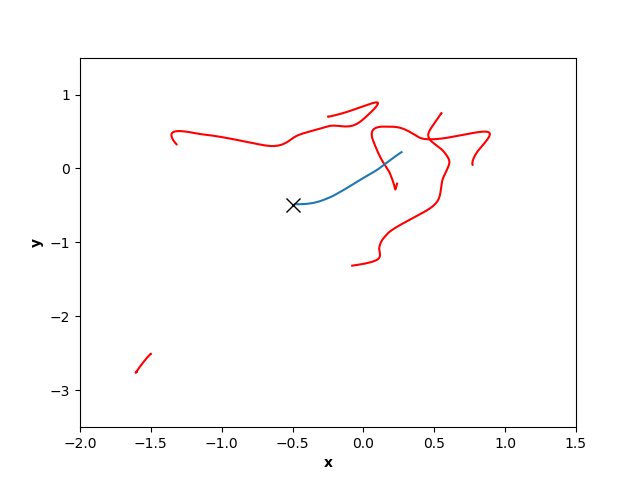}\\
	\includegraphics[scale=0.38]{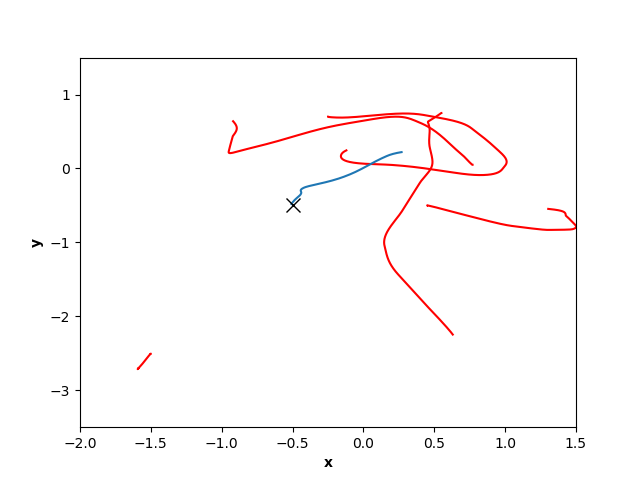}
	\includegraphics[scale=0.38]{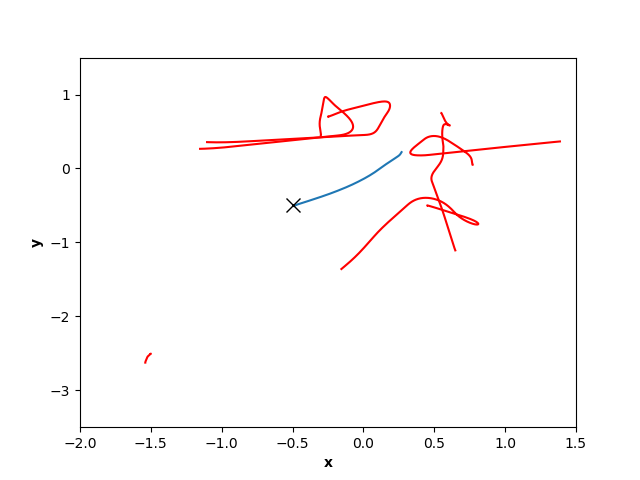}
	\caption{Space mapping trajectories for $1,\dots,6$ dogs. The simulations are stopped if $| \bar X - Z_\text{des}| < 0.05$ holds true. The corresponding times are $T_1 =3400, T_2 = 80, T_3 = 60, T_4 = 60, T_5 = 70, T_6 = 40,$ where the subscript refers to the number of dogs involved in the simulation. We see that already one dog is able to steer the crowd. Nevertheless, more dogs significantly decrease the time needed for the steering process. The stochastic influence in the system is implicitly displayed in the trajectory of the dogs in the figure on the top left. As for deterministic systems one would expect to have a homogeneous helix.}
	\label{fig:trajectories}
\end{figure}

Figure~\ref{fig:trajectories} compares the trajectories of the dogs (red) and the resulting trajectory of the centre of mass of the crowd (blue). Note,  that one dog has a hard time of leading the crowd as the iteration stops at $T_1 = 3400.$ The situation is getting better for two dogs. They are successful at time $T_2 = 80.$ Three dogs finish at time $T_3 = 60.$ In the other cases we have  $T_4 = 60, T_5 = 70, T_6 = 40.$ 

\begin{rem}
We emphasize that the initial positions of the dogs were chosen manually and not included in the optimization. Hence, we cannot deduce the optimal number of dogs from these results. 
\end{rem}

\subsection{Stabilization}
\begin{figure}[ht!]
	\includegraphics[scale=0.4]{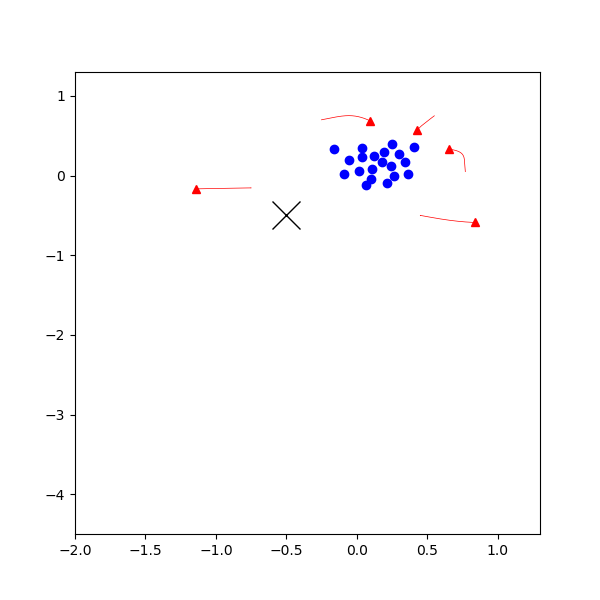}
	\includegraphics[scale=0.4]{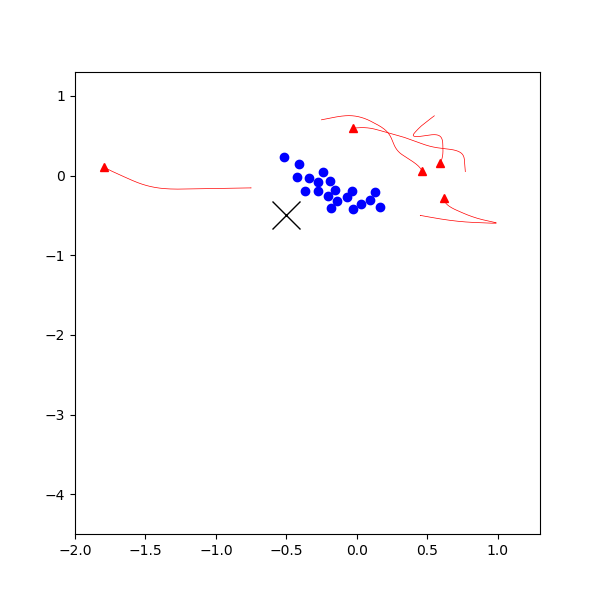}\\
	\includegraphics[scale=0.4]{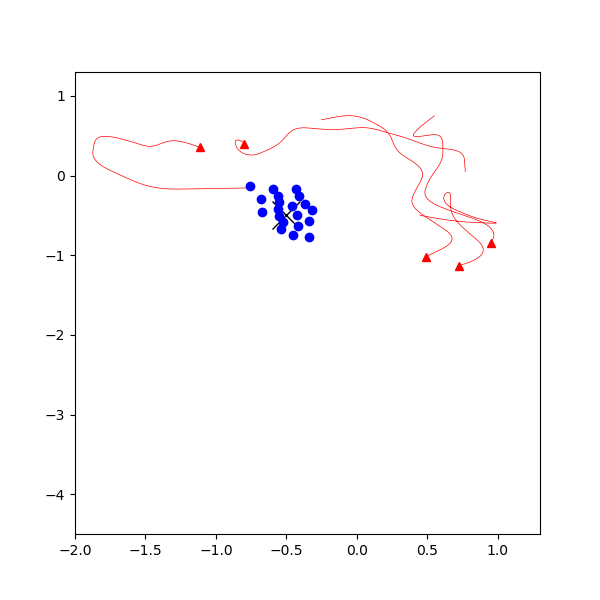}
	\includegraphics[scale=0.4]{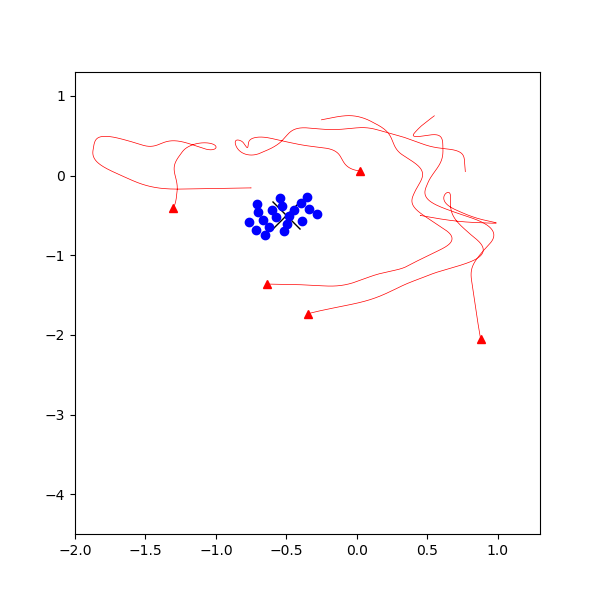}\\
	\includegraphics[scale=0.4]{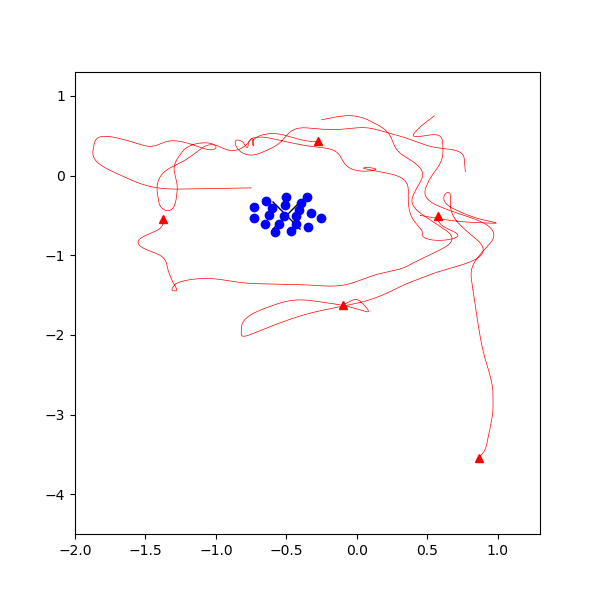} 
	\includegraphics[scale=0.4]{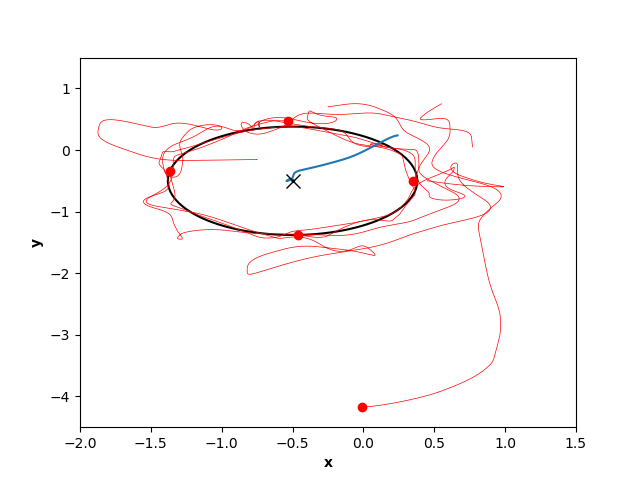}
	\caption{Snapshots of the optimization procedure at $t=10, 25, 50, 75, 125, 250$ (top-left to bottom-right). We see that the five dogs are able to steer the centre of mass of the crowd to the destination and that the crowd stays together. The latter is a new information which is not accessible by investigating only the centre of mass. Moreover, we see a stabilization as the dogs begin to circle around the crowd ($t=75, 125, 250$) after the centre of mass reached the desired destination. The black circle is underlining the discussion of the stabilization of the crowd.}
	\label{fig:simulation}
\end{figure}

Next, we show snapshots of a simulation with $T = 250, \gamma = 1e^{-3}$ and $5$ dogs in order to investigate if the herding process stabilizes. Indeed, we see in Figure~\ref{fig:simulation} that the dogs begin to circle around the crowd when the task of steering the centre of mass to the destination $Z_\text{des}$  is achieved. This behaviour can be interpreted as stabilization of the system. For this simulation the maximum number of space mapping iterations was limited to two. The stabilization indicated by this example is underlined by the following computations. For simplicity we set the friction parameter to $\alpha=0$ and consider the deterministic model.
	
\textbf{Herding one sheep with two dogs.}	
We first investigate the case having two dogs and only one sheep. Suppose the sheep is located at the destination $Z_\text{des}$ and the dogs are initially positioned at a circle with radius $\alpha p_1$ around $Z_\text{des}$, one at each end of a diameter. Then, the positions of the dogs can be parametrized with the help of $p_1$ as
$$ a_1 = Z_\text{des} + \alpha p_1, \qquad a_2 = Z_\text{des} - \alpha p_1.$$
In this setting the deterministic state equations are given by
\begin{align*}
\dot{x}_1 &= v_1, \qquad \dot{v}_1 = -G_2(|x_1 - a_1|) - G_2(|x_1-a_2|), \\
\dot{a_i} &= u_i, \quad i=1,2.
\end{align*}
Choosing the initial conditions $x_1 = Z_\text{des}, v_1(0) = 0, a_1 = Z_\text{des} + \alpha p_1, a_2 = Z_\text{des} - \alpha p_1,$ we obtain due to the radial symmetry of the interaction potentials
$$ \dot{v}_1 = -G_2(|x_1 - a_1|) - G_2(|x_1-a_2|) = 0. $$
Thus, the system is stable for $u_i = 0.$ Note, that we have reduced the system from two controls to only one control affecting both dogs at the same time. One can even change the positions of the dogs by $\dot{p}_1 = f$ for some function $f$ without any effect on the position of the sheep. This shows the stability of the configuration in this toy example.

\textbf{Herding $2N$ sheep with $2M$ dogs.} The observation of the previous section can be generalized to the following framework having an even number of dogs $2M$ and an even number of sheep $2N$ involved.

We assume the initial configuration to fulfil the following assumptions:
	(A1) $v_i = 0, \quad i=1,\dots,2N.$ \\
	(A2) The centre of mass is located at the destination $Z_\text{des}$, i.e.
	\[
	\frac{1}{2N} \sum_{i=1}^{2N} x_i = Z_\text{des}.
	\]
	(A3) For each $x_i$ exists $x_j$ with $i \ne j$ and $x_i - a_k = -(x_j - a_\ell)$ where $a_k$ and $a_\ell$ are positioned at a circle around $Z_\text{des}$ each on one end of a diameter.

Further, we note that we have assumed the interaction potentials to be radially symmetric. Using these assumptions we find
\begin{align*}
\frac{d}{dt}&\left(\frac{1}{2N} \sum_{i=1}^{2N} x_i \right) = \frac{1}{N} \sum_{i=1}^{2N} v_i \\ &= \frac{1}{N} \sum_{i=1}^{2N} \left( v_i (0) + \int_0^t \frac{1}{N} \sum_{j=1}^{2N} G_1(|x_i - x_j|) + \sum_{k=1}^{2M} G_2(|x_i - a_k|) \; \mathrm{d} s\right) \\
&= 0.
\end{align*}

Thus, setting $u_k =0$ for all $k=1,\dots,2M$ we obtain stability for the location of the centre of mass. Moreover, changing the velocities of the dogs as discussed in the framework with one sheep and two dogs, the position of the centre of mass is conserved as well. These findings are illustrated by the black circle added to the terminal configuration of the simulation shown in Figure~\ref{fig:simulation}. The dog at the bottom of the picture is far away from the crowd such that its contribution to the forces can be neglected. 

Of course, the stochastic influence and the error in the numerical integration will lead to dogs circling around the crowd. Thus, we do not expect to obtain numerically a stable setting with all sheep and dogs standing still.

\section{Space-mapping using the mean-field limit}\label{sec:meanfield}
Crowds consisting of many individuals, i.e.~$N\gg1$, are often investigated from a mesoscopic point of view with the help of a mean-field equation, see e.g.~\cite{AlbiPareschi,jose12,Schafe2,viewcones1}. Following the steps in \cite{paperMPC}, this equation can be obtained via the empirical density
\[
f^N(t,x,v) = \frac{1}{N}\sum_{i=1}^N \delta_0(x - x_i(t))\otimes \delta_0(v - v_i(t)).
\]
Formally passing to the mean-field limit $N\to \infty$ leads the deterministic optimization problem
\begin{equation}\label{eq:costfunctional-mf}
J(f,u ; \bar Z,\bar u) :=  \int_0^T \frac{1}{2N} \sum_{k=1}^N \| \E [f(t)] - \bar Z(t)\|^2  + \frac{\gamma}{2} \|u(t)- \bar u(t)\|_{\R^{MD}}^2 \dd t, 
\end{equation}
subject to the state system given by
\begin{align*}
\partial_t f + v\cdot \nabla_x(f) &= \div_v\Big( (G_1 \ast f) + \sum_{k=1}^M G_2(x - a_m) + \alpha v )\, f \Big), \\
\frac{d}{dt} a_m &= u_m, \qquad m=1,\dots,M.
\end{align*}
with initial conditions $f(0,x,v) = f_0(x,v),\; a_m(0) = a_0^m$ for $m=1,\dots,M.$
Here we used
$$  \E [f(t)] = \int_\mathbb{R} x\,f(t,x,v) \mathrm{d}x\,\mathrm{d}v.$$
Now, drawing the random initial conditions i.i.d.~from $f_0(x,v),$ it is well-known that $f(t,x,v)$ assigns the probability of finding a particle at time $t$ at position $x$ with velocity $v.$ Hence, one could use the deterministic mean-field problem for $f$ as coarse model for a space-mapping in order to control the stochastic limit for many particles. Note that a similar deterministic optimization problem was solved in \cite{paperMPC}. 

\begin{rem}
We want to emphasize that in the common noise case the limiting equation is not deterministic but a stochastic PDE \cite{Lacker}. Thus, it is not clear whether it is an appropriate choice as coarse model. The cost for computing optimal controls with the SPDE are probably very high.
\end{rem}

\section{Conclusion and Outlook} \label{sec:ConOutlook}
We discussed a space mapping approach in combination with receding horizon control for the closed loop control of a stochastic interacting particle system. The numerical results underline that the method is feasible for interacting particle systems with small stochastic perturbation. Further, they indicate that a sub-optimal control for the stochastic system is found efficiently already after few space mapping iterations. 

In near future, we plan to use the space mapping approach to control a stochastic system involving a large number of interacting particles. In this case,  the mean-field approximation can be used as coarse model for the space mapping approach. Moreover, a rigorous analysis of the space mapping procedure applied to stochastic problems is of interest.

Further, an investigation of the stabilizing effect of the feedback control is planned. A rigorous generalization to the setting with common noise and its influence on the space-mapping performance are interesting future projects as well.

\bibliographystyle{spmpsci}      
\bibliography{biblio}
\end{document}